\documentclass[12pt,a4paper]{amsart}
\usepackage{amsfonts}
\usepackage{amsthm}
\usepackage{amsmath}
\usepackage{amscd}
\usepackage[latin2]{inputenc}
\usepackage{t1enc}
\usepackage[mathscr]{eucal}
\usepackage{indentfirst}
\usepackage{graphicx}
\usepackage{graphics}
\usepackage{pict2e}
\usepackage[T1]{fontenc}
\usepackage{currvita}
\usepackage{url}
\usepackage{epic}
\numberwithin{equation}{section}
\usepackage[margin=2.9cm]{geometry}
\usepackage{epstopdf}

\theoremstyle{plain}

 \theoremstyle{definition}

\newtheorem{?}[Th]{Problem}

\begin{document}
\author{Vojt\u{e}ch Dvo\u{r}\'ak}
\address{Trinity College, Cambridge CB21TQ, UK. Email address: vd273@cam.ac.uk}

\title{A note on Norine's antipodal-colouring conjecture
}


\begin{abstract} 
Norine's antipodal-colouring conjecture, in a form given by Feder and Subi,
asserts that whenever the edges of the discrete cube are 2-coloured there must exist a
path between two opposite vertices along which there is at most one colour change.
The best bound to date was that there must exist such a path with at most $n/2$ colour
changes. Our aim in this note is to improve this upper bound to $(\frac{3}{8}+o(1))n$. 
\end{abstract}

\maketitle

\section{Introduction} 

The hypercube $ Q_{n} $ has vertex set $   \lbrace 0,1 \rbrace ^{n} $, with two vertices joined by an edge if they differ in a single coordinate. We call two vertices of $Q_{n}$ antipodal if their graph distance is $n$. We call a pair of edges of $Q_{n}$ $v_{1}w_{1}$ and $v_{2}w_{2}$ antipodal if $v_{1}$ and $v_{2}$ are antipodal vertices and $w_{1}$ and $w_{2}$ are also antipodal vertices. A 2-colouring of the edges of $ Q_{n} $ is called \textit{antipodal} if no pair of antipodal edges has the same colour. Norine \cite{norine} conjectured the following.

\vspace{5mm}
\textbf{Conjecture 1.} In any antipodal 2-colouring of the edges of $ Q_{n} $, there
exists a pair of antipodal vertices which are joined by a monochromatic path.
\vspace{5mm}

Feder and Subi \cite{feder} later made the following conjecture.

\vspace{5mm}
\textbf{Conjecture 2.} In any 2-colouring of the edges of $ Q_{n} $, we can find a pair of antipodal vertices and a path joining them with at most one colour change.

\vspace{5mm}
If true, this implies the conjecture of Norine. Indeed, consider an antipodal 2-colouring of the edges of $ Q_{n} $. By Conjecture 2, we can now find an antipodal path $ P_{1}P_{2} $ such that both paths $ P_{1} $ and $ P_{2} $ are monochromatic. If they have the same colour we are done; if not the path $ P_{1}P_{2}^{C} $ will work, where $ P_{2}^{C} $ is the antipodal path to path $ P_{2} $.

Call a path in $ Q_{n} $ a \textit{geodesic} if no two of its edges have the same direction. Leader and Long \cite{leader} proved the following result.

\vspace{5mm}
\textbf{Theorem 3.} In any 2-colouring of the edges of $ Q_{n} $, we can find a monochromatic geodesic of the length of at least $ {\lceil}\frac{n}{2}{\rceil} $.
\vspace{5mm}

Leader and Long proposed a conjecture that strengthens that of Feder and Subi.

\vspace{5mm}
\textbf{Conjecture 4.} In any 2-colouring of the edges of $ Q_{n} $, we can find a pair of antipodal vertices and a geodesic joining them with at most one colour change.

\vspace{5mm}

Theorem 3 implies we can always find a pair of antipodal vertices and a geodesic joining them with at most $ \frac{n}{2} $ colour changes. Moreover, as Theorem 3 is sharp, there is no hope of improving the result by finding longer monochromatic geodesic. In the present note, we establish.

\vspace{5mm}
\textbf{Theorem 5.} In any 2-colouring of the edges of $ Q_{n} $, we can find a pair of antipodal vertices and a geodesic joining them with at most $ (\frac{3}{8}+o(1))n $ colour changes.
\vspace{5mm}

To do that, we employ the strategy of dividing the $Q_{n}$ graph into small pieces ($ Q_{3} $ graphs in fact) and finding collection of geodesics of certain properties within each piece. The conditions we impose on these local geodesics let us glue them together into collection of geodesics of $Q_{n}$ in such a way that on average these long geodesics will have not too many colour changes. From that we in particular conclude that at least one of the long geodesics must have not too many colour changes. In the next section we collect together some simple facts about $Q_{3}$ graphs, and then we use these in the following section to prove Theorem 5.

\section{Good and bad $ Q_{3} $ graphs}

In this section we collect together some facts about 2-colourings of the 3-dimensional cube. Of course, each one is trivial, but we record them here as they will be useful in what follows.

From now on, we call a geodesic connecting two antipodal points simply an \textit{antipodal geodesic}.

We call a colouring of $ Q_{3} $ by two colours \textit{good} if we can find 4 antipodal geodesics, with each vertex being an endpoint of exactly one of these, such that these 4 geodesics have in total at most two colour changes. If a colouring of $ Q_{3} $ is not good we call it \textit{bad}.

The terms good and bad $ Q_{3} $ will be sometimes used instead of good and bad colouring of $ Q_{3} $, and it is understood that we refer to a particular colouring.

When showing that Conjecture 4 holds for n=5, Feder and Subi \cite{feder} proved the following simple lemma which we will use too.

\vspace{5mm}
\textbf{Lemma 6. } Assume in a 2-colouring of $ Q_{3} $, there are antipodal vertices $ v $ and $ v' $ such that all the geodesics connecting $ v $ and $ v' $ have two colour changes. Then the other three pairs of antipodal points are connected by geodesics without colour changes.
\vspace{5mm}

In particular, note that this implies that in any bad $ Q_{3} $, any pair of antipodal vertices is connected by a geodesic with at most one colour change. So it follows that in bad $ Q_{3} $, at most one pair of antipodal points is connected by a geodesic without colour change.

Our first lemma gives us an easy way to identify many cube graphs as good.

\vspace{5mm}
\textbf{Lemma 7. }Suppose all three edges at some vertex of $ Q_{3} $ have the same colour. Then it is a good colouring.

\begin{proof}
Assume this colouring is bad. Without loss of generality take the vertex where all edges have the same colour to be 000 and this colour to be blue. If all the edges with neither of their endpoints being 000 or 111 are red, it is a good colouring, as the other three pairs of antipodal vertices are connected by the antipodal geodesics with no colour changes. So assume some edge with neither endpoint being 000 or 111 is blue, without loss of generality it is (100,110).

From 001, we have the antipodal geodesic with no colour change (001,000,100,110). So if it is bad, for no other pair of antipodal points we can have an antipodal geodesic with no colour change. So the edge (100,101) must be red by considering the geodesic (010,000,100,101), the edge (001,101) must be red by considering the geodesic (010,000,001,101) and the edge (001,011) must be red by considering the geodesic (100,000,001,011).

But that gives the red antipodal geodesic with no colour change (100,101,001,011), thus a contradiction. 
\end{proof}

Next, note that one particular example of a bad cube graph occurs when we colour all the edges in one direction by one colour, and all the edges in the other two directions by the other colour. Lemma 8 that follows tells us that any bad colouring behaves very much like this example in the sense we will need in our proof.

\vspace{5mm}
\textbf{Lemma 8. }Consider any bad colouring of $ Q_{3} $ and any vertex $ v $. Then there exists an antipodal geodesic from $ v $ to $ v' $ with exactly one colour change, red edge at $ v $ and blue edge at $ v' $.
\begin{proof}
Without loss of generality let $ v $ be 000. Assume no such antipodal geodesic exists. By Lemma 7, we have at least one red edge from 000, without loss of generality to 100. Also we have at least one blue edge from 111. If this edge went to 110 or 101, we would be done immediately, so it must go to 011 and the other two edges from 111 must be red. Furthermore, the other two edges from 000 must be blue, else we would be done, so assume they are blue.

As we are in a bad $ Q_{3} $, at most one pair of antipodal points can be connected by a antipodal geodesic without a colour change. Whichever colour the edge (001,101) has, it creates an antipodal geodesic without a colour change, either between 010 and 101 or between 001 and 110. So 000 cannot be connected to 111 by an antipodal geodesic without a colour change, forcing the edges (100,110) and (100,101) to be blue and the edges (010,011) and (001,011) to be red.

But now we see that both 010 and 001 are connected to their antipodals by geodesics without a colour change, contradiction.
\end{proof}

Analogously, in any bad $ Q_{3} $, there exists such an antipodal geodesic with exactly one colour change, blue edge at $ v $ and red edge at $ v' $.

\section{Proof of the main result}

We will now prove the main result.

\vspace{5mm}
\textbf{Theorem 5.} In any 2-colouring of edges of $ Q_{n} $, we can find a pair of antipodal vertices and a geodesic joining them with at most $ (\frac{3}{8}+o(1))n $ colour changes.

\begin{proof}
As we have $ o(n) $ term included in our bound, clearly it suffices to prove the theorem for $ n $ divisible by 3, so assume $ n=3k $. For the vertices $ v,w $ of distance 3, call $ G(v,w) $ the subgraph of $ Q_{n} $ spanned by the geodesics between $ v $ and $ w $ (so $ G(v,w) \cong  Q_{3} $). Call two such subgraphs $ G_{1} \cong Q_{3} $ and $ G_{2} \cong Q_{3} $ of $ Q_{n} $ \textit{neighbours} if they share exactly one vertex. If this vertex is $ v $, call them  \textit{v-neighbours}. Consider a set $ A $ of all the ordered pairs $ (v,w) $ of the vertices of $ Q_{n} $ such that $ d(v,w)=3 $. Assume $ f:A \rightarrow V(Q_{n}) $ satisfies the following three conditions for all the vertices $ v, w $:

(i)$d(v,f(v,w))=1  $

(ii)$d(w,f(v,w))=2  $

(iii)$d(f(w,v),f(v,w))=1  $

In other words, this is equivalent to $ (v,f(v,w),f(w,v),w) $ being an antipodal geodesic in $ G(v,w) $.

Now, given the antipodal geodesic $ (v_{0},v_{1},v_{2},v_{3},v_{4},...,v_{3i},v_{3i+1},v_{3i+2},v_{3i+3},...,v_{3k}) $, we will modify it into the antipodal geodesic $ (v_{0},f(v_{0},v_{3}),f(v_{3},v_{0}),v_{3},f(v_{3},v_{6}),
...,v_{3i},\\f(v_{3i},v_{3i+3}),f(v_{3i+3},v_{3i}),v_{3i+3},...,v_{3k}) $. 

I will show that for every fixed colouring, we can define such $ f $ in a way that the expected number of colour changes on a random modified antipodal geodesic is no more than $ (\frac{3}{8}+o(1))n $, where $ o(1) $ term depends on $ n $ only, not on the colouring. More precisely, I will define $ f_{1} $ and $ f_{2} $ (depending on the colouring) and show that at least one of these two must always work.

If $ G(v,w) $ is a good $ Q_{3} $, for $ i=1,2 $, let $ f_{i}(v,w) $ and $ f_{i}(w,v) $ be same for both values of $ i $ and such that no other geodesic between $ v $ and $ w $ has strictly less colour changes than $ (v,f_{i}(v,w),f_{i}(w,v),w) $.

Call the vertex $ v $ of $ Q_{n} $ \textit{even} if its distance from 000...000 is even and call it \textit{odd} otherwise. Every geodesic of the length 3 connects an odd and an even vertex. For $ G(v,w) $ bad with $ v $ even and $ w $ odd,  define $ f_{1}(v,w) $ and $ f_{1}(w,v) $ such that $ (v,f_{1}(v,w),f_{1}(w,v),w) $ has exactly one colour change, $ (v,f_{1}(v,w)) $ is blue and $ (f_{1}(w,v),w) $ is red. Also define $ f_{2}(v,w) $ and $ f_{2}(w,v) $ such that $ (v,f_{2}(v,w),f_{2}(w,v),w) $ has exactly one colour change, $ (v,f_{2}(v,w)) $ is red and $ (f_{2}(w,v),w) $ is blue. By Lemma 8, there exist such functions.

Denote by $p$ the proportion of good $ Q_{3} $ subgraphs of $ Q_{n} $ in this colouring and the proportion of bad ones is thus $ 1-p $. Picking two $ Q_{3} $ subgraphs that are neighbours uniformly at random, denote the probability that both are good by $ a $, the probability that one is good and one is bad by $ b $, and thus the probability that both are bad is $ 1-a-b $. We clearly must have $ p=a+\frac{b}{2} $. 

How large can $ b $ be? Suppose at any vertex $ v $, of all $ Q_{3} $ containing $ v $, there is $ s $ good ones and $ {n\choose 3} - s $ bad ones. There are $ \frac{1}{2} {n\choose 3} {n-3\choose 3} $ pairs of v-neighbours, and of them at most $ s({n\choose 3}-s) \leq \frac{1}{4} {n\choose 3}^{2} $ are good-bad pairs. As this applies to every vertex and is independent of $ s $, we have $ b \leq \frac{1}{2} {n\choose 3}{n-3\choose 3}^{-1}=\frac{1}{2}+o(1) $.

Now, choose a modified antipodal geodesic uniformly at random. Due to the symmetry, and the properties of good and bad cube graphs, for any $ j: 0\leq j \leq n-1 $ and for either value of $i$, the expected number of colour changes inside the geodesic $ (v_{3j},f_{i}(v_{3j},v_{3j+3}),f_{i}(v_{3j+3},v_{3j}),v_{3j+3}) $ is at most $ \frac{1}{2}p+(1-p)=1-\frac{p}{2} $.

What is the probability that, for some fixed $ j: 1\leq j \leq n-1 $, we have a colour change between the edges $ (f_{i}(v_{3j},v_{3j-3}),v_{3j}) $ and  $ (v_{3j},f_{i}(v_{3j},v_{3j+3})) $ (due to the symmetry this is same for all such $ j $)? With probability $ 1-a-b $, both $ G(v_{3j},v_{3j-3}) $ and $ G(v_{3j},v_{3j+3}) $ are bad, and then we do not have a colour change by definition of $ f_{i} $. If one is good and one is bad, exactly one of $ f_{1} $ and $ f_{2} $ has a colour change between these two edges. So choose as our $ f $ that $ f_{i} $ for which the probability of a change in this case is at most $ \frac{1}{2} $. 

Finally, with probability $ a $, both graphs are good. Consider any fixed vertex $ v $. Choosing a random subgraph $ Q_{3} $ containing $ v $, by Lemma 7, the probability that it is good is at least the probability that choosing 3 random distinct edges from $ v $, they all have the same colour. So we conclude there are at least $ (\frac{1}{4}-o(1))n^{3} \geq \frac{1}{8}n^{3} $ good subgraphs containing $ v $ for $ n $ large enough. Suppose precisely $ t $ good subgraphs contain $ v $. Clearly, the number of pairs of neighbours of good subgraphs that have a colour change at $ v $ is at most $ \frac{1}{4}t^{2} $. Also for any good subgraph $ G_{1} $ containing $ v $, the number of good subgraphs that share $ v $ and at least one other vertex with $ G_{1} $ is less than $ 3n^{2} $. So the number of pairs of two good graphs that are v-neighbours is at least $ \frac{1}{2}t(t-3n^{2}) $. So the probability that a uniform random pair of good v-neighbours switches colour there is at most $ \frac{1}{2}\frac{t}{t-3n^{2}} $. This is a decreasing function of $ t $ for $ t>3n^{2} $, so using $ t \geq \frac{1}{8}n^{3} $, and as this applies to any vertex, we get that in this case, the probability of a colour switch is no more than $ \frac{1}{2}+o(1) $.

Putting these together we obtain that on average our modified antipodal geodesic has at most $ (1-\frac{p}{2})\frac{n}{3}+(\frac{b}{2}+(\frac{1}{2}+o(1))a)\frac{n}{3} =(\frac{1}{3}+\frac{b}{12}+o(1))n$ colour changes (using $ p=a+\frac{b}{2} $). But $ b \leq \frac{1}{2} + o(1) $, giving the result.
\end{proof}

\section*{Acknowledgements}

I would like to thank Imre Leader for suggesting this problem and for various helpful discussions, my PhD supervisor B\'ela Bollob\'as for advice regarding the final version of this note and to Trinity College for supporting the research.


\begin{thebibliography}{9}


\bibitem{feder}
T. Feder and C. Subi.
\textit{On hypercube labellings and antipodal monochromatic paths}.
Discrete Appl. Math. \textbf{161} (2013), 1421-1426.

\bibitem{leader}
I. Leader and E. Long.
\textit{Long geodesics in subgraphs of the cube}.
Discrete Math. \textbf{326} (2014), 29-33.


\bibitem{norine} 
S. Norine. 
\textit{Edge-antipodal colorings of cubes}.
Open Problem Garden, \url{http://www.openproblemgarden.org/op/edge_antipodal_colorings_of_cubes}.





\end{thebibliography}
\end{document}